\documentclass{article}
\usepackage[utf8]{inputenc}
\pdfoutput=1

\usepackage[margin=1.25in]{geometry}

\usepackage{amsmath,color,graphicx,amssymb}
\usepackage{amsthm}
\usepackage{amsmath}
\usepackage{tikz}
\usetikzlibrary{calc,matrix,arrows,shapes.geometric,positioning,decorations.pathmorphing}
\usepackage{tikz-cd}

\usepackage{float}
\usepackage{slashbox}
\usepackage{indentfirst}
\usepackage{fancyhdr}
\usepackage{enumitem}

\fancypagestyle{firstpage}{\fancyhf{}\fancyfoot[L]{\textit{Keywords}: Koszul homology, Koszul algebra, edge ideals}}

\newtheorem{thm}{Theorem}[section]

\theoremstyle{definition}
\newtheorem{defn}[thm]{Definition}

\theoremstyle{definition}
\newtheorem{ex}[thm]{Example}

\theoremstyle{definition}

\theoremstyle{definition}

\theoremstyle{definition}
\newtheorem{quest}[thm]{Question}

\theoremstyle{definition}
\newtheorem{conj}[thm]{Conjecture}

\allowdisplaybreaks

\title{\large{\uppercase{A survey on the Koszul homology algebra}}}
\author{\normalsize{\uppercase{Rachel N. Diethorn}}}
\date{}

\begin{document}

\maketitle

\begin{abstract}
The Koszul homology algebra of a commutative local (or graded) ring $R$ tends to reflect important information about the ring $R$ and its properties.  In fact, certain classes of rings are characterized by the algebra structure on their Koszul homologies.  In this paper we survey some classical results on the Koszul homology algebras of such rings and highlight some applications.  We report on recent progress on the Koszul homology algebras of Koszul algebras and examine some open questions on the topic.
\end{abstract}

\section{Introduction}

The properties of a local (or graded) ring $R$ are often encoded in the Koszul homology algebra $H(R)$; and in certain cases, the key to understanding the ring is to understand the algebra structure on its Koszul homology.  Roughly speaking, the simpler the algebra structure on $H(R)$, the nicer the ring $R$.  For example, a local ring $R$ with residue field $k$ is regular if and only if $H(R)=k$; this follows from the characterization of regular local rings given by Eilenberg and by Auslander and Buchsbaum (see for example \cite[Theorem 1.4.13]{gulliksenlevin} or \cite[Lemma 5]{tate}).  Complete intersection rings, Gorenstein rings, and Golod rings are also characterized by their Koszul homology algebras; we discuss these classical characterizations in Section 4.  Such characterizations have been quite useful in a wide range of applications in commutative algebra, particularly so in the study of Golod rings, and very recently in calculating Poincar{\'e} series over certain classes of rings; see for example \cite{rossisega} and \cite{ksv}.  We highlight a few such applications at the end of Section 5. 

The classical characterizations discussed in this survey have also played important roles in the development of complete classifications of the Koszul homology algebras of quotients of regular local rings by certain ideals.  Such classifications were given by Avramov, Kustin, and Miller in \cite{akm} and Weyman in \cite{weyman} for ideals of grade 3, by Kustin and Miller in \cite{kustingor} for Gorenstein ideals of grade 4, and by Kustin in \cite{kustin} for almost complete intersections of grade 4.  We discuss these results and a few of their applications in Section 5.

While there are complete classifications of the algebra structures on Koszul homologies in certain cases, for a general ring, the structure can be quite complicated.  This is evidenced by the paper \cite{roos}, in which Roos determines completely the Poincar{\'e} series (and the graded Betti numbers) of 104 rings with embedding dimension 4 whose defining ideals are generated by quadratic forms, and examines the different homological properties of their Koszul homology algebras.  Even for Koszul\footnote{We use the word ``Koszul" in two different ways throughout this survey: \textit{Koszul homology algebras} $H(R)$, which we define in Section 2, and \textit{Koszul algebras} $R$ which we define in Section 6.} algebras, whose Koszul homology algebras are known to be restricted in various ways, some of which we describe in Section 6, Roos shows that the structure of $H(R)$ can be surprisingly complicated.  We highlight one such result in Theorem \ref{thmroos} and one remarkable example in Example \ref{crolletalex1}, and leave the rest of \cite{roos} for the reader to explore.

Although Koszul algebras are not known to be characterized completely by the algebra structures on their Koszul homologies, the strong connections between the Koszul property of $R$ and the Koszul homology algebra of $R$ has provided an active and vibrant area of research in recent years.  The results of and questions posed by Avramov, Conca, and Iyengar in their papers \cite{aci1} and \cite{aci2} were catalysts for much of this progress, and thus are the focus of Section 6 of this survey.  In particular we discuss three questions about the Koszul homology algebras of Koszul algebras; the first is a question of Avramov about the generation of $H(R)$ as an algebra; the second and third are questions of Avramov, Conca, and Iyengar about bounds on the Betti numbers of $R$, and about subadditivity of syzygies of $R$, respectively.  We observe that the latter two questions have equivalent formulations in terms of Koszul homology, and we outline the progress made towards answering these questions and what remains unanswered.  

A recurring theme throughout this survey is the use of the Koszul homology algebra $H(R)$ as a tool for understanding the ring $R$; of course this requires having effective methods for studying Koszul homology.  One such approach, which has been used in the proofs of several important results that we state in this paper, is that of minimal models, which is developed by Avramov in \cite[Section 7.2]{6lectures}.  Another approach, which was first developed by Herzog in \cite{herzog92}, is to find, explicitly, the generators of Koszul homology modules.  One of the key ingredients of Herzog's approach is an isomorphism that relates Koszul homologies $H_i(R)$ to the minimal resolution of $R$.  We discuss this isomorphism and state Herzog's result in Section 3, and we refer back to this machinery often throughout the paper.

\section{Preliminaries}

In this section we recall the definition of the Koszul homology algebra of a graded ring (the definition in the local case is similar).  Let $R=\bigoplus_{i\geq 0}R_i$ be a standard graded $k$-algebra with $k$ a field, and fix a minimal set of generators $\underline{x}=x_1,\dots,x_n$ of $R_1$.  In this survey we focus our attention on the \textit{Koszul complex} $K(\underline{x})$, which is the complex with modules 
\begin{align*}
K_i=\bigwedge\nolimits^i(R^n)
\end{align*}
for $0\leq i\leq n$ and $K_i=0$ for all $i>n$, where we fix a basis $dx_1,\dots,dx_n$ of the free module $R^n$, and differentials
\begin{align*}
\partial^K_i(dx_{j_1}\cdots dx_{j_i})=\sum_{\ell=1}^i(-1)^{\ell+1}x_{j_{\ell}}dx_{j_1}\cdots\widehat{dx_{j_{\ell}}}\cdots dx_{j_i},
\end{align*} 
where here we shorten $K(\underline{x})$ and $K(\underline{x})_i$ to $K$ and $K_i$, respectively.  We note that $K_0=R$, $K_1=R^n$, and for $2\leq i\leq n$, $K_i$ is a free $R$-module of rank $n\choose i$ with basis
\begin{align*}
\{dx_{j_1}\cdots dx_{j_i}\,|\, 1\leq j_1<\dots<j_i\leq n\}.
\end{align*} 

Given our fixed minimal set of generators $\underline{x}$ of $R_1$, we define $K(R)$ and $H_i(R)$ to be the Koszul complex $K(\underline{x})$ and its $i^{\mathrm{th}}$ homology $H_i(K(\underline{x}))$, respectively; although, it is worth noting that for any other minimal set of generators $\underline{y}$ of $R_1$, the resulting Koszul homology modules $H_i(K(\underline{y}))$ are isomorphic to $H_i(K(\underline{x}))$ (see for example \cite[Remark 1.4]{huneke}).  Recall that Koszul homology satisfies \textit{depth sensitivity}; that is, the last nonvanishing Koszul homology module is $H_{n-g}(R)$ where $g=\mathrm{depth}\,R$.  For more background on the Koszul complex, see for example \cite{matsumura}.  

We denote by $H(R)$, the sum of the Koszul homology modules

\begin{align*}
H(R)=\underset{i\geq 0}{\bigoplus} H_i(R).
\end{align*}
The Koszul complex is a \textit{differential graded algebra}, or a \textit{DG algebra}; that is, it is an algebra whose differential satisfies the Leibniz rule:
\begin{align*}
\partial^K(ab)=\partial^K(a)b+(-1)^{|a|}a\partial^K(b), 
\end{align*} 
where $|a|$ denotes the homological degree of $a$.  One can easily verify this fact using the definition given above.  The DG algebra structure on $K(R)$ induces a $k$-algebra structure on $H(R)$.  Indeed, given two cycles $a$ and $b$ in $K(R)$, the equality $\partial^K(ab)=0$ follows directly from the Leibniz rule; thus, $ab$ is a cycle, and one can define an algebra structure on $H(R)$ by
\begin{align}\label{alg}
[a]\cdot[b]:=[ab], 
\end{align} 
where $[\,\cdot\,]$ denotes the homology class.  From (\ref{alg}) it is easy to see that $H(R)$ inherits \textit{graded commutativity}; that is, $y\cdot z=(-1)^{|y||z|}z\cdot y$ for all $y$ and $z$ in $H(R)$, from $K(R)$.  We call $H(R)$ the \textit{Koszul homology algebra} of $R$ (even though it is defined only up to non-canonical isomorphism of $k$-algebras).     

\section{Canonical Bases}

Let $k$ be a field and $R$ a standard graded $k$-algebra.  Then $R\cong Q/I$ where $Q=k[x_1,\dots,x_n]$ is a standard graded polynomial ring with $I$ a homogeneous ideal.  The Koszul homology algebra $H(R)$ is closely connected to the the minimal free resolution $F$ of $R$ over $Q$ via the isomorphsim
\begin{align}\label{eq:1}
H_i(R)=H_i(R\underset{Q}{\otimes}K(Q))=\mathrm{Tor}_i^Q(R,k)=H_i(F\underset{Q}{\otimes}k)\cong F_i\underset{Q}{\otimes}k.  
\end{align}
As such, Koszul homology has been used classically as a tool for studying the minimal free resolution of $R$ and vice versa.  We refer to this isomorphism often throughout this survey.  

In \cite{herzog92}, Herzog uses the isomorphism above to explicitly describe canonical $k$-bases of the Koszul homology modules, in the case that $k$ is a field of characteristic zero.  The basis elements are described in terms of Jacobians that arise from the differentials in the minimal graded free resolution $F$ of $R$, as we see in the following result.

To set notation for the theorem, for each $i$ we let $b_i$ be the rank of $F_i$,  and fix a homogeneous basis $e_1^i,\dots,e_{b_i}^i$ of $F_i$ such that $\partial^F(e_j^i)=\sum_{\ell=1}^{b_{i-1}} f_{\ell,j}^i e_{\ell}^{i-1}$.  For $z\in K(Q)$, we denote by $\overline{z}$ the image of $z$ in $R\otimes K(Q)$.  In the theorem, $c_{j_1,\dots,j_i}$ is the inverse of a constant that tracks the degrees of the corresponding basis elements of $F$.    

\begin{thm}\label{herzoggens}(\cite[Corollary 2]{herzog92})   
For each $i=1,\dots,n$, a $k$-basis of $H_i(R)$ is given by $[\overline{z}_{j}]$ for $j=1,\dots,b_i$, where 
\begin{align*}
z_{j}=\sum_{1\leq \ell_1<\dots<\ell_i\leq n}\sum_{j_1=1}^{b_1}\dots\sum_{j_{i-1}=1}^{b_{i-1}}c_{j_1,\dots,j_i}\frac{\partial\left(f_{j_{i-1},j}^i,f_{j_{i-2},j_{i-1}}^{i-1},\dots,f_{j_1,j_2}^2,f_{1,j_1}^1\right)}{\partial\left(x_{\ell_1},x_{\ell_2},\dots,x_{\ell_i}\right)}dx_{\ell_1}\dots dx_{\ell_i}.
\end{align*}
\end{thm}

This result of Herzog was generalized in \cite{CGHPU} and \cite{DIETHORN2020106387}; however, the description of the basis elements given in both generalizations require that $k$ is a field of characteristic zero.  In \cite{Herzog2018}, Herzog and Maleki give a different description of a basis for $H_i(R)$ that does not depend on the characteristic of the field.  This description is given in terms of certain operators on $Q$, which are defined as follows.  For $f\in(x_1,\dots,x_n)$ and for $r=1,\dots,n$, let
\begin{align*}
d^r(f)=\frac{f(0,\dots,0,x_r,\dots,x_n)-f(0,\dots,0,x_{r+1},\dots,x_n)}{x_r}.
\end{align*}
It is clear that these operators are $k$-linear maps and that they depend on the order of the variables.  

\begin{thm}\label{thmhm}(\cite[Theorem 1.3]{Herzog2018})
For each $i=1,\dots,n$, a $k$-basis of $H_i(R)$ is given by $[\bar{z}_{j}]$ for $j=1,\dots,b_i$, where 
\begin{align*}
z_{j}=\sum_{1\leq \ell_1<\dots<\ell_i\leq n}\sum_{j_1=1}^{b_1}\dots\sum_{j_{i-1}=1}^{b_{i-1}}d^{\ell_i}(f_{j_{i-1},j}^i)d^{\ell_{i-1}}(f_{j_{i-2},j_{i-1}}^{i-1})\dots d^{\ell_2}(f_{j_1,j_2}^2)d^{\ell_1}(f_{1,j_1}^1)dx_{\ell_1}\dots dx_{\ell_i}.
\end{align*}
\end{thm} 

It is important to note that both theorems above give generators for the modules $H_i(R)$, but do not directly provide information about the algebra structure on $H(R)$.  However, such bases have been used as tools to understand the Koszul homology algebra; see for example Theorem \ref{diethornedgeideals}, where Theorem \ref{thmhm} is used to study the Koszul homology algebras of certain Koszul algebras.  The theorems above have also been particularly useful in obtaining criteria for Golodness, which we discuss in Section 4.2 of this survey. 

We finish this section by observing that just as the isomorphism (\ref{eq:1}) can be used to describe Koszul homologies $H_i(R)$ in terms of the data in the minimal free resolution of $R$, as was done in Theorems \ref{herzoggens} and \ref{thmhm}, it can also be used to describe the minimal resolution of $R$ in terms of its Koszul homologies.  The latter was the focus of \cite{aramovaherzog95} and \cite{aramovaherzog96}, in which Aramova and Herzog describe explicitly how to construct the minimal free resolution of $R$, once bases of Koszul homologies are known.  These constructions provide our first concrete example of how the Koszul homologies of $R$ can be used as a tool for learning about the ring $R$, and we will see several more examples throughout this survey. 

\section{Classical Characterizations}

In this section we discuss classical characterizations of complete intersection rings, Gorenstein rings, and Golod rings by their Koszul homology algebras.  Such characterizations have played a tremendous role in the effort to understand these classes of rings.  We highlight a few direct applications of these characterizations throughout the section.   

\subsection{Complete Intersection and Gorenstein Rings}

In this section we state the classical results that characterize complete intersection and Gorenstein rings in terms of their Koszul homology algebras.  We begin with the characterization of complete intersection rings due to Assmus and Tate; the equivalence of the first three conditions is due to Assmus, and the equivalence of the first and last is due to Tate.  

\begin{thm}(\cite[Theorem 2.7]{assmus},\,\,\cite[Theorem 6]{tate})\label{CI}
The following conditions are equivalent: 
\begin{description}
\item[(1)]  $R$ is a local complete intersection;
\item[(2)]  $H(R)$ is isomorphic to the exterior algebra on $H_1(R)$;
\item[(3)]  $H_2(R)=(H_1(R))^2$;
\item[(4)]  $H(R)$ is generated by $H_1(R)$ as a $k$-algebra.
\end{description}
\end{thm}

Soto generalizes Theorem \ref{CI} slightly in \cite[Proposition 3]{soto}; the result asserts that $R$ is a complete intersection if and only if $H(R)$ is a free graded $k$-algebra.  Another characterization of complete intersection rings by their Koszul homology algebras is given by Bruns.  In order to state this characterization we first recall that there is a natural map
\begin{align*}
\lambda\colon\bigwedge H_1(R)\longrightarrow H(R)
\end{align*}
which extends the identity map on $H_1(R)$; this fact follows from the universal property of the exterior algebra on $H_1(R)$, or from the graded commutativity of $H(R)$.  Thus it follows from Theorem \ref{CI} that $R$ is a complete intersection ring if and only if the map $\lambda$ is surjective.  The result of Bruns below describes complete intersections by the injectivity of this map.

\begin{thm}(\cite[Theorem 2]{bruns})\label{bruns}
If $R$ is a local ring containing a field, then the following statements are true: 
\begin{description}[style=unboxed,leftmargin=0cm]
\item[(1)]  $H_1(R)^i=0$ for $i>\mathrm{edim}\,R-\mathrm{dim}\,R$;
\item[(2)]  $R$ is a complete intersection if and only if $\lambda$ is injective.
\end{description}
\end{thm}

This result is applied in the proof of Theorem \ref{bhici}, which provides a condition under which the defining ideal of a Koszul algebra is a complete intersection.  

An analagous characterization of Gorenstein rings in terms of their Koszul homology algebras involves the condition known as Ponicar{\'e} duality, which we define below.

\begin{defn}
A graded algebra $H=\bigoplus_{i=0}^n H_i$ satisfies \textit{Poincar{\'e} duality} if for each $i=0,\dots,n$, the $H_0$-homomorphism
\begin{align*}
H_i&\longrightarrow\mathrm{Hom}_{H_0}(H_{n-i},H_n) \\
h&\longmapsto\phi_h
\end{align*} 
with $\phi_h(x)=hx$, is an isomorphism.
\end{defn}

Now we state the result of Avramov and Golod which asserts that Gorenstein rings are characterized by their Koszul homology algebras.

\begin{thm}(\cite{avramovgolod})\label{gor}  Let $R$ be a local ring, and let $n=\mathrm{edim}\,R-\mathrm{depth}\,R$.  The following conditions are equivalent:  
\begin{description}[style=unboxed,leftmargin=0cm]
\item[(1)]  $R$ is a Gorenstein ring;
\item[(2)]  $H(R)$ satisfies Poincar{\'e} duality;
\item[(3)]  the $k$-linear map $H_{n-1}(R)\longrightarrow\mathrm{Hom}_k(H_1(R),H_n(R))$ induced by the multiplication on $H(R)$ is injective.
\end{description}
\end{thm}

We highlight applications of Theorem \ref{gor} at the end of Section 4.2 and in Section 5.

\subsection{Golod Rings}

In this section we discuss a classical result which characterizes Golod rings by their Koszul homology algebras, and a more recent result which stems from this classical theorem.  At the end of the section, we highlight an application of Theorem \ref{gor} which provides a connection between the Golod and Gorenstein properties.  We begin with a definition of a Golod ring.  

\begin{defn}
A local ring $R$ is a \textit{Golod ring} if its Poincar{\'e} series satisfies the equality
\begin{align}\label{eq:2}
P_k^R(t)=\frac{(1+t)^{\mathrm{edim}\,R}}{1-\sum_{j=1}^n\mathrm{rank}_k H_j(R)t^{j+1}},
\end{align}
where $n=\mathrm{edim}\,R-\mathrm{depth}\,R$. 
\end{defn}

A classical result of Serre asserts that, for any local ring, the formula in (\ref{eq:2}) provides a coefficient-wise upper bound for its Poincar{\'e} series; see for example \cite{golod}, \cite[Theorem 1.3]{levin1976lectures}, \cite[Proposition 3.3.2]{MR2641236}, or \cite[Corollary 4.2.4]{gulliksenlevin}.  Thus, roughly speaking, the resolution of the residue field $k$ of a Golod local ring has the fastest growth possible.   

In \cite{golod}, Golod introduced some higher order operations, called Massey operations, on Koszul homology, and used them to characterize rings satisfying the equality (\ref{eq:2}), which we now call Golod rings.  For a more detailed treatment of Massey operations, see for example \cite[Section 4.2]{gulliksenlevin}.

\begin{defn}
We say that $K(R)$ admits a \textit{trivial Massey operation} if for some $k$-basis $\mathcal{B}=\{h_{\lambda}\}_{\lambda\in\Lambda}$ of $H_{\geq 1}(R)$, there is a function
\begin{align*}
\mu:\bigsqcup_{n=1}^{\infty}\mathcal{B}^n\rightarrow K(R) 
\end{align*}
such that $\mu(h_{\lambda})=z_{\lambda}$ is a cycle with $[z_{\lambda}]=h_{\lambda}$ and
\begin{align}\label{eq:3}
\partial^K\mu(h_{\lambda_1},\dots,h_{\lambda_p})=\sum_{j=1}^{p-1}\overline{\mu(h_{\lambda_1},\dots,h_{\lambda_j})}\mu(h_{\lambda_{j+1}},\dots,h_{\lambda_p})
\end{align}
where $\overline{a}=(-1)^{|a|+1}a$.
\end{defn}

We state Golod's classical result using the modern terminology of Golod rings, as follows.

\begin{thm}(\cite{golod})\label{golod}
A local ring $R$ is Golod if and only if $K(R)$ admits a trivial Massey operation.
\end{thm}

Thus Golod rings are characterized by having highly trivial Koszul homology algebras.  Indeed, if $K(R)$ admits a trivial Massey operation, then (\ref{eq:3}) implies that the product on $H(R)$ is trivial; that is, $H_{\geq 1}(R)\cdot H_{\geq 1}(R)=0$.  Note, however, that a trivial product on $H(R)$ is not enough to imply that $R$ is Golod; the higher Massey products must also be trivial.  In \cite{katthan}, Katth{\"a}n produces an example which illustrates this fact.  

\begin{thm}(\cite[Theorem 3.1]{katthan})
Let $Q=k[x_1,x_2,y_1,y_2,z]$ be a polynomial ring with $k$ a field, let $I$ be the ideal
\begin{align*}
I=(x_1x_2^2,\,\, y_1y_2^2,\,\, z^3,\,\, x_1x_2y_1y_2,\,\, y_2^2z^2,\,\, x_2^2z^2,\,\, x_1y_1z,\,\, x_2^2y_2^2z),
\end{align*}
and let $R=Q/I$.  Then the product on $H(R)$ is trivial, but $R$ is not Golod.
\end{thm}  

For the remainder of this section, we recognize a few direct applications of the characterizations above.  First, we look at an application of Theorem \ref{golod} and Theorem \ref{herzoggens} to Golod rings.  In \cite{herzoghuneke} Herzog and Huneke exploit the canonical bases of Koszul homologies given in Theorem \ref{herzoggens} and the characterization of Golod rings in Theorem \ref{golod} to provide a differential condition for Golodness, and in doing so, they are able to produce large classes of Golod rings, including quotients by powers and symbolic powers of ideals.  The differential condition is given in the following result; see also \cite[Sections 2 \& 3]{Herzog2018}, \cite[Theorem 3.5]{gupta}, and \cite[Proposition 4.4]{DIETHORN2020106387} for similar applications. 

\begin{thm}
Let $Q$ be a standard graded polynomial ring over a field of characteristic zero, let $I\subseteq Q$ be a homogeneous ideal, and let $R=Q/I$.  Let $\partial(I)$ denote the ideal generated by the partial derivatives of the elements of $I$.  If $(\partial(I))^2\subseteq I$, then $R$ is a Golod ring. 
\end{thm}

Now we look at an application of Theorem \ref{gor} by Avramov and Levin in \cite{avramovlevin} which establishes a connection between Golod homorphisms and Gorenstein rings.  The following theorem asserts that factoring a Gorenstein ring by its socle produces a Golod homomorphism; that is, a relative version of the Golod property.  We discuss an application of Theorem \ref{avramovlevin} to calculating Poincar{\'e} series in the following section.

\begin{thm}(\cite[Theorem 2]{avramovlevin})\label{avramovlevin}
Let $R$ be a local Gorenstein ring with $\mathrm{edim}\,R>1$ and $\mathrm{dim}\,R=0$.  Then
\begin{align*}
R\longrightarrow R/(0:m)
\end{align*}
is a Golod homomorphism.
\end{thm}  

Throughout this section, we have seen several classes of local rings which are completely characterized by the algebra structures on their Koszul homology algebras, and we conclude by emphasizing that, even for an arbitrary local ring, the Koszul homology algebra is such a sensitive invariant that it encodes all of the numerical information about the free resolution of its residue field, as we see in the following result of Avramov.

\begin{thm}(\cite[Corollary 5.10]{hopf})
The Poincar{\'e} series $P_k^R(t)$ of a local ring $R$ with residue field $k$ depends only on the algebra structure on $H(R)$ and its higher order Massey operations.
\end{thm}

\section{Classifications of Koszul homology algebras}

For a general ring, the algebra structure on its Koszul homology can be quite complicated (see for example \cite{roos}); however, there are certain classes of rings whose Koszul homology algebras are completely classified.  We outline those cases in this section.

We begin with the case of $R=Q/I$, where $Q$ is a local ring with residue field $k$, whose defining ideal $I$ is perfect of grade 3; that is, the length of the longest $Q$-regular sequence contained in $I$ and the projective dimension of $R$ over $Q$ are both equal to 3.  In this case, Weyman \cite{weyman} and Avramov, Kustin, and Miller \cite{akm} (see also \cite{avramovclassify}) provide complete classifications of the algebra $\mathrm{Tor}_{\bullet}^Q(R,k)$.  Restricting to the case where $Q$ is regular, this gives the following classification of the possible multiplicative structures on the Koszul homology algebra $H(R)$ via the isomorphism (\ref{eq:1}).  Notice that in this case the last nonvanishing Koszul homology module is $H_3(R)$ by the Auslander-Buchsbaum formula and depth sensitivity of the Koszul complex.

\begin{thm}(\cite[Theorem 2.12]{akm})\label{akm}
Let $R=Q/I$ with $Q$ a regular local ring.  If the projective dimension of $R$ over $Q$ is 3, then there are nonnegative integers p, q, and r and bases $\{e_i\}$, $\{f_i\}$, and $\{g_i\}$ of $H_1(R)$, $H_2(R)$, and $H_3(R)$, respectively, such that the multiplication on $H(R)$ is given by one of the following:
\begin{description}[style=unboxed,leftmargin=0cm]
\item[CI]:  $f_1=e_2e_3,\,\, f_2=e_3e_1,\,\, f_3=e_1e_2,\,\, e_if_j=\delta_{ij}g_1$\, for $1\leq i,j\leq 3$
\item[TE]:  $f_1=e_2e_3,\,\, f_2=e_3e_1,\,\, f_3=e_1e_2$
\item[B]:  $e_1e_2=f_3,\,\, e_1f_1=g_1,\,\, e_2f_2=g_1$
\item[G($r$)]:  $e_if_i=g_1$\, for $1\leq i\leq r$\, and $r\geq 2$
\item[H($p,q$)]:  $e_{p+1}e_i=f_i$\, for $1\leq i\leq p$\, and $e_{p+1}f_{p+i}=g_i$\, for $1\leq i\leq q$
\end{description}
with $e_je_i=-e_ie_j$, $e_i^2=0$, and $e_if_j=f_je_i$ for all $i$ and $j$, and with the products of basis elements that are not listed above being zero.
\end{thm}  

The class listed as \textbf{CI} is precisely the class where $I$ is a complete intersection; the class listed as \textbf{G($r$)} is the class where $I$ is Gorenstein, but not a complete intersection, and $r$ indicates the minimal number of generators of $I$; the class listed as \textbf{TE} consists of rings that are neither Golod, nor complete intersections.

The classification in Theorem \ref{akm} builds on the classification for perfect grade 3 almost complete intersection ideals given by Buchsbaum and Eisenbud in \cite[Theorem 5.3]{BE}; see also \cite{brown} for another special case of this classification.  

Such classifications for the Tor algebra are also known for grade 4 Gorenstein ideals  and grade 4 almost complete intersection ideals; the former is due to Kustin and Miller \cite[Theorem 2.2]{kustingor}, and the latter is due to Kustin \cite[Theorem 1.5]{kustin}.  Again, restricting to the case where the ambient ring is regular, we present the former case here, and we note that the proof applies the characterization of Gorenstein rings in Theorem \ref{gor}.    

\begin{thm}(\cite[Theorem 2.2]{kustingor})\label{kustingor}
Let $R=Q/I$ with $Q$ a regular local ring, and assume that every element in $k$ has a square root in $k$.  If $I$ is a grade 4 Gorenstein ideal that is not a complete intersection, then there are bases $e_1,\dots,e_n$ for $H_1(R)$; $f_1,\dots,f_{n-1},f_1',\dots,f_{n-1}'$ for $H_2(R)$; $g_1,\dots,g_n$ for $H_3(R)$; and $h$ for $H_4(R)$ such that the multiplication $H_i(R)\cdot H_{4-i}(R)$ is given by
\begin{align*}
e_ig_j = \delta_{ij}h,\,\, f_if_j'=\delta_{ij}h,\,\, f_if_j=f_i'f_j'=0,
\end{align*}  
with all other products given by one of the following cases:
\begin{description}[style=unboxed,leftmargin=0cm]
\item[(1)]  $H_1(R)\cdot H_1(R)=0$ and $H_1(R)\cdot H_2(R)=0$
\item[(2)]  All products in $H_1(R)\cdot H_1(R)$ and $H_1(R)\cdot H_2(R)$ are zero except:
\begin{displaymath}
e_1e_2=f_3,\,\, e_1e_3=-f_2,\,\, e_2e_3=f_1,\,\, e_1f_2'=-e_2f_1'=g_3,\,\, -e_1f_3'=e_3f_1'=g_2,\,\,e_2f_3'=-e_3f_2'=g_1
\end{displaymath}
\item[(3)]  There is an integer $p$ such that $e_{p+1}e_i=f_i$, $e_if_i'=g_{p+1}$, and $e_{p+1}f_i'=-g_i$ for $1\leq i\leq p$, and all other products in $H_1(R)\cdot H_1(R)$ and $H_1(R)\cdot H_2(R)$ are zero. 
\end{description}
\end{thm}

Not surprisingly, we see that the algebra structure in the grade 4 case is more complicated than that in the grade 3 case.  

The classifications of Koszul homology algebras in Theorems \ref{akm} and \ref{kustingor} have been instrumental in the study of rationality of Poincar{\'e} series, which has been a topic of intense interest over the years and remains a central problem in homological commutative ring theory today; see for example \cite[Section 6]{mp} or \cite[Section 4.3]{MR2641236} for a complete treatment of this problem.  For example, Theorem \ref{akm} is an important ingredient in the proof of \cite[Theorem 6.4]{akm}, which determines several classes of rings over which the Poincar{\'e} series of a finitely generated module is rational; analagously, \cite[Theorem A]{jkm} applies the classification in Theorem \ref{kustingor}.  In fact, much of the recent progress on rationality of Poincar{\'e} series relies on the classical characterizations of Koszul homology discussed in Section 4.  As a direct application of Theorem \ref{avramovlevin}, Rossi and {\c{S}}ega calculate the Poincar{\'e} series of a class of artinian Gorenstein rings in \cite[Proposition 6.2]{rossisega} and show that they are indeed rational.  For similar applications of the Koszul homology algebra structure in calculating Poincar{\'e} series, see \cite[Theorem 5.1]{rossisega} and \cite[Theorem 7.1]{ksv}.  

We finish this section by recognizing one more recent application of the classifications above.  In \cite{christensenveliche}, Christensen and Veliche determine minimal cases for which powers of the maximal ideal in a local ring are not Golod.  We state their result in embedding dimension 3, which uses the classification in Theorem \ref{akm}; the embedding dimension 4 case \cite[Proposition 5.2]{christensenveliche} uses Theorem \ref{kustingor}.

\begin{thm}(\cite[Theorem 4.2]{christensenveliche})
Let $R$ be an artinian Gorenstein local ring of embedding dimension 3 and socle degree 3 with maximal ideal $\mathfrak{m}$.  The following conditions are equivalent:
\begin{description}[style=unboxed,leftmargin=0cm]
\item[(1)]  $R$ is a complete intersection;
\item[(2)]  $R$ is compressed and Koszul;
\item[(3)]  $R/\mathfrak{m}^3$ belongs to the class \textbf{TE};
\item[(4)]  $R/\mathfrak{m}^3$ is not Golod.
\end{description}
\end{thm} 

Such applications demonstrate the usefulness of the Koszul homology algebra $H(R)$ as a tool for learning about $R$, and discovering its properties.

\section{Recent Progess: Koszul Algebras}

In this section, we focus on recent progress on the Koszul homology algebras of \textit{Koszul algebras}; that is, $k$-algebras over which $k$ has a linear resolution.  One can check that the defining ideals of Koszul algebras are generated by quadratics; however, not all algebras defined by quadratics are Koszul (see for example \cite[Remark 1.10]{concasurvey}).  Classical examples of Koszul algebras include quotients of polynomial rings by quadratic monomial ideals (e.g. edge ideals), Veronese algebras, and Segre product algebras.  Unlike the classes of rings discussed in Section 4, Koszul algebras are not known to be characterized solely by the algebra structures on their Koszul homologies; nonetheless, the Koszul property of $R$ is closely connected to the algebra structure on its Koszul homology $H(R)$.

To present these connections, we assume throughout this section that $R=Q/I$ is a standard graded $k$-algebra with $Q=k[x_1,\dots,x_n]$ and $k$ a field (although some results can be stated more generally), and we view the Koszul homology algebra $H(R)$ of a Koszul algebra $R$ as a bigraded algebra
\begin{align*}
H(R)=\underset{i,j}{\bigoplus}H_i(R)_j,
\end{align*}   
where $i$ is the homological degree, and $j$ is the internal degree given by the grading on $R$.  Given the isomorphism (\ref{eq:1}), the rank of $H_i(R)_j$ is given by the graded Betti number $\beta_{ij}$ of $R$.  As such, we can view these bigraded pieces of the Koszul homology algebra in the following table, which comes from the Betti table of $R$:
\begin{center}
\begin{tabular}{ |c|c c c c c  } 
 \hline
 \backslashbox{$j-i$}{$i$} & 0 & 1 & 2 & 3 & \dots\\ 
 \hline
 0 & $H_{0,0}$ & $H_{1,1}$ & $H_{2,2}$ & $H_{3,3}$ & \dots \\
 1 & $H_{0,1}$ & $H_{1,2}$ & $H_{2,3}$ & $H_{3,4}$ & \dots \\ 
 2 & $H_{0,2}$ & $H_{1,3}$ & $H_{2,4}$ & $H_{3,5}$ & \dots \\ 
 3 & $H_{0,3}$ & $H_{1,4}$ & $H_{2,5}$ & $H_{3,6}$ & \dots \\ 
 \vdots & \vdots & \vdots & \vdots & \vdots & 
\end{tabular}
\end{center}
where $H_{i,j}=H_i(R)_j$.  

It is customary to call row one of the Betti table the \textit{linear strand} and the other rows \textit{nonlinear strands}.  More precisely, we have the following definition.

\begin{defn}
We call the subspace
\begin{align*}
\underset{i}{\bigoplus}\,H_i(R)_{i+1}
\end{align*}
of $H(R)$ the \textit{linear strand} of $H(R)$.
\end{defn}

Now we recall some useful terminology for describing the shapes of Betti tables, which we will use throughout this section.  In \cite{aci1}, Avramov, Conca, and Iyengar introduce the sequence
\begin{align*}
t_i(R)=\mathrm{sup}\{j\in\mathbb{Z}\,|\,H_i(R)_j\neq 0\},
\end{align*}
which encodes important information about $R$.  For example, the integer $t_i(R)-i$ is the height of the $i$th column of the Betti table.  Notice that the \textit{regularity} of $R$
\begin{align*}
\mathrm{reg}(R)=\underset{i\geq 0}{\mathrm{sup}}\{t_i(R)-i\}
\end{align*}
measures the height of the Betti table above.  The \textit{slope} of $R$,
\begin{align*}
\mathrm{slope}(R)=\underset{i\geq 1}{\mathrm{sup}}\left\{\frac{t_i(R)-t_0(R)}{i}\right\},
\end{align*}
introduced and studied in \cite{aci1}, measures the slope of the Betti table.  The last number we will use in this section to describe the shapes of Betti tables is 
\begin{align*}
m(R)=\mathrm{min}\{i\in\mathbb{Z}\,|\,t_i(R)\geq t_{i+1}(R)\},
\end{align*} 
which is introduced in $\cite{aci2}$.  

The first result we state is a special case of a more general result of Avramov, Conca and Iyengar in \cite[Main Theorem]{aci1}; the inequality in (2) follows from results of Backelin in \cite{backelin} and of Kempf in \cite[Lemma 4]{kempf}. 

\begin{thm}\label{aci1mainthm} 
If $R$ is Koszul, then $R$ satisfies the following conditions:
\begin{description}[style=unboxed,leftmargin=0cm]
\item[(1)]  $\mathrm{slope}(R)=2$;
\item[(2)]  $t_i(R)\leq 2i$ for all $i\geq 0$;
\item[(3)]  $\mathrm{reg}(R)\leq\mathrm{pd}_QR$.
\end{description}
\end{thm}

The general result from which Theorem \ref{aci1mainthm} follows is one about the slope and regularity of $R$ where $Q$ is Koszul, with no assumptions on $R$, and the proof is based on the theory of minimal models; see for example \cite[Section 7.2]{6lectures}.  We will see several other applications of minimal models throughout this section.  

It follows from Theorem \ref{aci1mainthm} that if $R$ is Koszul, then $H_i(R)j=0$ for all $j>2i$, which we record in Theorem \ref{aci2mainthm}, and one can easily check that $H_i(R)_i=0$ for all $i\geq 1$; in other words, the Betti table above has the form:
\begin{center}
\begin{tabular}{ |c|c c c c c  } 
 \hline
 \backslashbox{$j-i$}{$i$} & 0 & 1 & 2 & 3 & \dots\\ 
 \hline
 0 & $H_{0,0}$ & $0$ & $0$ & $0$ & \dots \\
 1 & $0$ & $H_{1,2}$ & $H_{2,3}$ & $H_{3,4}$ & \dots \\ 
 2 & $0$ & $0$ & $H_{2,4}$ & $H_{3,5}$ & \dots \\ 
 3 & $0$ & $0$ & $0$ & $H_{3,6}$ & \dots \\ 
 \vdots & \vdots & \vdots & \vdots & \vdots & 
\end{tabular}
\end{center}

Theorem \ref{aci1mainthm}, along with the other results of \cite{aci1} and \cite{aci2}, sparked an intensive investigation of the Koszul homology algebra of a Koszul algebra, which we outline throughout the rest of this section.

\subsection{Generation by the linear strand}

In this section we discuss a question of Avramov about the generation of the Koszul homology algebra of a Koszul algebra.  We begin with the following result of Avramov, Conca, and Iyengar which records a consequence of Theorem  \ref{aci1mainthm} and describes the algebra structure on the main diagonal of the Betti table.  This result prompted Question \ref{quest}, which we discuss throughout this section.   

\begin{thm}(\cite[Theorem 5.1]{aci2})\label{aci2mainthm}
If $R$ is Koszul, then 
\begin{description}
\item[(1)]  $H_i(R)j=0$ for all $j>2i$, and
\item[(2)]  $H_i(R)_{2i}=(H_1(R)_2)^i$ for all $i\geq 0$.
\end{description}
\end{thm}

In fact, Avramov, Conca, and Iyengar prove a more general version of Theorem \ref{aci2mainthm} where the Koszul hypothesis is relaxed to only require linearity in the resolution of $k$ up to degree $n$.  Under this relaxed hypothesis, the equalities in the theorem hold for all $i\leq n-1$.  Again, the proof is based on the theory of minimal models.

By Theorem \ref{aci2mainthm}, one can see that the elements of the Koszul homology algebra which lie on the main diagonal are contained in the subalgebra generated by the linear strand.  This observation led Avramov to ask the following question.

\begin{quest}\label{quest}
If $R$ is Koszul, is the Koszul homology algebra of $R$ generated as a $k$-algebra by the linear strand?  
\end{quest}

In \cite[Theorem 3.1]{bdgms}, the authors extend Theorem \ref{aci2mainthm} to the next diagonal of the Betti table; their result is stated as follows.

\begin{thm}(\cite[Theorem 3.1]{bdgms})\label{bdgmsgen}
If $R$ is Koszul, then 
\begin{align*}
H_i(R)_{2i-1}=(H_1(R)_2)^{i-2}H_2(R)_3
\end{align*}
for all $i\geq 2$.
\end{thm}

However, the answer to Question \ref{quest} is negative in general.  The first counterexample was found computationally by Eisenbud and Caviglia on Macaulay2 \cite{M2}, and this example led Conca and Iyengar to consider quotients by edge ideals of $n$-cycles.  From this, Boocher, D'Al{\'i}, Grifo, Monta{\~n}o, and Sammartano were able to produce the family of counterexamples in (3) of the following theorem.  The proofs of both Theorems \ref{bdgmsgen} and \ref{bdgms} continue the trend of using the minimal model of $R$ over $Q$ to study $H(R)$. 

\begin{thm}(\cite[Theorem 3.15]{bdgms})\label{bdgms}
Let $Q=k[x_1,\dots,x_n]$, let $I$ be the edge ideal associated to a graph $G$ on the vertices $x_1,\dots,x_n$, and let $R=Q/I$.
\begin{description}[style=unboxed,leftmargin=0cm]
\item[(1)]  If $G$ is an $n$-path, then $H(R)$ is generated by $H_1(R)_2$ and $H_2(R)_3$.
\item[(2)]  If $G$ is an $n$-cycle with $n\not\equiv 1(\mathrm{mod}\,3)$, then $H(R)$ is generated by $H_1(R)_2$ and $H_2(R)_3$. 
\item[(3)]  If $G$ is an $n$-cycle with $n\equiv 1(\mathrm{mod}\,3)$, then for any $0\neq z\in H_{\lceil\frac{2n}{3}\rceil}(R)_n$, $H(R)$ is generated by $H_1(R)_2$, $H_2(R)_3$, and $z$. 
\end{description}
\end{thm}

In \cite[Theorem 4.2]{diethornedgeideals} the author applies Theorem \ref{thmhm}, and extends  Theorem \ref{bdgms}(1) to give a broader class of edge ideals for which  Question \ref{quest} has a positive answer.

\begin{thm}(\cite[Theorem 4.2]{diethornedgeideals})\label{diethornedgeideals}
Let $Q=k[x_1,\dots,x_n]$ with $k$ a field, let $I$ be the edge ideal associated to a forest on the vertices $x_1,\dots,x_n$, and let $R=Q/I$.  Then $H(R)$ is generated by the linear strand. 
\end{thm}

The following result of Conca, Katth{\"a}n, and Reiner gives another class of Koszul algebras whose Koszul homology algebras are generated by the linear strand.

\begin{thm}(\cite[Corollary 2.4]{concakatthanreiner})
Let $Q=k[x_1,\dots,x_n]=\bigoplus_i Q_i$ with $k$ a field of characteristic zero, and let $R=\bigoplus_i Q_{2i}$ be the second Veronese subalgebra.  Then $H(R)$ is generated by the linear strand.
\end{thm}

As the authors of \cite{concakatthanreiner} point out, Question \ref{quest} remains open for other classes of Koszul algebras, such as higher Veronese subalgebras and Segre product algebras.  Furthermore, the authors of \cite{bdgms} ask the following weaker version of Question \ref{quest}.

\begin{quest}(\cite[Question 3.16]{bdgms})\label{questbdgms}
Is there a Koszul algebra $R$ which is a domain and whose Koszul
homology $H(R)$ is not generated as a $k$-algebra by the linear strand?
\end{quest}

The answer to this question is affirmative; the example below was communicated to the authors of \cite{bdgms} by McCullough.

\begin{ex}
Let $k$ be a field and let
\begin{align*}
R=k[ae,\, af,\, ag,\, ah,\, bg,\, bh,\, ce,\, cg,\, ch,\, de,\, df,\, dg].
\end{align*}
Notice that $R$ is the toric edge ring of a bipartite graph, hence a domain, and its defining ideal is generated by quadratic binomials; indeed, $R\cong k[X_1,X_2,X_3,X_4,X_5,X_6,X_7,X_8,X_9,X_{10},X_{11},X_{12}]/I$, where
\begin{align*}
I=\,&\big(X_3X_{11}-X_2X_{12},\, X_8X_{10}-X_7X_{12},\, X_3X_{10}-X_1X_{12},\, X_2X_{10}-X_1X_{11}, \\
&X_6X_8-X_5X_9,\, X_4X_8-X_3X_9,\, X_4X_7-X_1X_9,\, X_3X_7-X_1X_8,\, X_4X_5-X_3X_6\big),
\end{align*} 
thus $R$ is Koszul by \cite{oh}.  Using Macaulay2 \cite{M2}, we see that the Betti table of $R$ is 
\begin{center}
\begin{tabular}{ |c|c c c c c c } 
 \hline
  & 0 & 1 & 2 & 3 & 4 & 5  \\ 
 \hline
 0 & $1$ & - & - & - & - & -  \\
 1 & - & $9$ & $11$ & - & - & -   \\ 
 2 & - & - & $10$ &  $26$ & $15$ & $1$ \\ 
 3 & - & - & - & - & - & $1$ \\  
\end{tabular}
\end{center}
and it is easy to see that the Koszul homology algebra must have a generator in bidegree $(5,7)$ for degree reasons.  Thus $H(R)$ is not generated by the linear strand.
\end{ex}

So far we have seen that Koszulness of $R$ does not imply generation by the linear strand in $H(R)$.  It is also the case that generation by the linear strand in $H(R)$ does not imply Koszulness of $R$, as we see in the example below.

\begin{ex}(\cite[7.4]{crolletal})\label{crolletalex1}
Let $k$ be a field, and let 
\begin{align*}
R=k[X,Y,Z,U]/(X^2 + XY, XZ + YU, XU, Y^2, Z^2, ZU + U^2)
\end{align*}
(see \cite[Case 55]{roos}).  $H(R)$ is generated by the linear strand; indeed, it has only 6 generators in bidegree $(1,2)$ and 4 generators in bidegree $(2,3)$.  However, $R$ is not Koszul; the resolution of $k$ is not linear.
\end{ex}

However, the authors of \cite{crolletal} provide stronger conditions on the algebra structure of $H(R)$ that are enough to imply Koszulness of $R$.  They prove that if $H(R)$ is generated by either a single element of bidegree $(1,2)$, or by a special set of elements in the linear strand, then $R$ is Koszul.  More precisely, their result is as follows.

\begin{thm}(\cite[Theorem 6.1]{crolletal})
Assume one of the following conditions holds:
\begin{description}[style=unboxed,leftmargin=0cm]
\item[(1)] There exists an element $[z]$ of bidegree $(1,2)$ such that every element in the nonlinear strands of $H(R)$ is a multiple of $[z]$.
\item[(2)]  $R_{\geq 3}=0$ and there is a set of cycles $Z$ representing elements in the linear strand with the property that $zz'=0$ for all $z,z'\in Z$ whose homology classes generate the nonlinear strands of $H(R)$.
\end{description}
Then R is Koszul.
\end{thm} 

In fact, they prove that the condition in (1) implies that $R$ is absolutely Koszul, a condition which implies Koszulness.

\subsection{Upper bounds on Betti numbers}

In this section we discuss the following question of Avramov, Conca, and Iyengar about upper bounds of the Betti numbers of $R$, and equivalently, by (\ref{eq:1}), the ranks of the Koszul homology modules of Koszul algebras.

\begin{quest}(\cite[Question 6.5]{aci1})\label{questbetti}
If $R=Q/I$ is Koszul, and $I$ is minimally generated by $g$ quadrics, does the inequality
\begin{align*}
\beta_i^Q(R)\leq {g\choose i},
\end{align*}
equivalently, $\mathrm{dim}\,_kH_i(R)\leq{g\choose i}$, hold for all $i$?  In particular, is the projective dimension of $R$ over $Q$ at most $g$?
\end{quest}

Standard arguments show that the answer to this question is positive for Koszul algebras whose defining ideals are generated by monomials or have Gr{\"o}bner bases of quadrics.  Furthermore, the answer to this question is positive for Koszul algebras whose defining ideals are minimally generated by $g\leq 4$ quadrics and for Koszul almost complete intersection algebras; that is, algebras with $\mathrm{ht}\,I=g-1$, as we outline below.  Otherwise, Question \ref{questbetti} remains open.  

The case where $g\leq 3$ is addressed by Boocher, Hassanzadeh, and Iyengar in the following result. 

\begin{thm}(\cite[Theorem 4.5]{boocherhiyengar})\label{boocherhiyengarbetti}
If $I$ is generated by 3 quadrics, then the following conditions are equivalent.
\begin{description}
\item[(1)]  $R$ is Koszul;
\item[(2)]  $H(R)$ is generated by the linear strand;
\item[(3)]  The Betti table of $R$ over $Q$ is one of the following:
\end{description}
\begin{center}
\begin{tabular}{ |c|c c c c  } 
 \hline
  & 0 & 1 & 2 & 3  \\ 
 \hline
 0 & $1$ & - & - & -  \\
 1 & - & $3$ & - & -   \\ 
 2 & - & - & $3$ &  - \\ 
 3 & - & - & - & $1$ \\  
\end{tabular}
\hspace{0.6cm}
\begin{tabular}{ |c|c c c c  } 
 \hline
  & 0 & 1 & 2 & 3  \\ 
 \hline
 0 & $1$ & - & - & -  \\
 1 & - & $3$ & $1$ & -   \\ 
 2 & - & - & $2$ &  $1$ \\ 
\end{tabular}
\hspace{0.6cm}
\begin{tabular}{ |c|c c c   } 
 \hline
  & 0 & 1 & 2   \\ 
 \hline
 0 & $1$ & - & -   \\
 1 & - & $3$ & $2$    \\ 
\end{tabular}
\hspace{0.6cm}
\begin{tabular}{ |c|c c c c  } 
 \hline
  & 0 & 1 & 2 & 3  \\ 
 \hline
 0 & $1$ & - & - & -  \\
 1 & - & $3$ & $3$ & $1$   \\ 
\end{tabular}
\end{center}
\end{thm}

In \cite{boocherhiyengar} the authors note that the first Betti table listed in Theorem \ref{boocherhiyengarbetti} corresponds to a complete intersection of three quadrics; the second is the Betti table of the ring $k[x,y,x]/(x^2,y^2,xz)$; the third corresponds to the ideal of minors of a $3\times 2$ matrix of linear forms (e.g. $k[x,y]/(x,y)^2$); the last is the Betti table of a ring with a linear resolution (e.g. $k[x,y,z]/(x^2,xy,xz)$).

The implication $(1)\implies (3)$ in Theorem \ref{boocherhiyengarbetti} gives an affirmative answer to Question \ref{questbetti} in the case where $g\leq 3$, and by the isomorphism (\ref{eq:1}), it determines the ranks of the bigraded pieces of $H(R)$ for such algebras; the implication $(1)\implies (2)$ provides another class of Koszul algebras for which Question \ref{quest} has a positive answer. 

The proofs of the implications $(2)\implies (1)$ and $(3)\implies (1)$ use D'Al{\'i}'s classification of Koszul algebras defined by 3 quadrics in \cite{dali}, which we state below.

\begin{thm}(\cite[Theorem 3.1]{dali})
Let $k$ be an algebraically closed field of characteristic different from 2, let $Q$ be a polynomial ring over $k$, and let $R=Q/I$ with $I$ a quadratic ideal of $Q$. If $\mathrm{dim}_kR_2=3$, then $R$ is Koszul if and only if it is not isomorphic as a graded $k$-algebra (up to trivial fiber extension) to any of these:
\begin{description}[style=unboxed,leftmargin=0cm]
\item[(1)]  $k[x,y,z]/(y^2+xy,\, xy+z^2,\, xz)$
\item[(2)]  $k[x,y,z]/(y^2,\, xy+z^2,\, xz)$
\item[(3)]  $k[x,y,z]/(y^2,\, xy+yz+z^2,\, xz)$.
\end{description}
In particular, if $R$ is a non-Koszul algebra with $\mathrm{edim}\,R=3$ defined by 3 quadrics, then its Betti table is
\begin{center}
 \begin{tabular}{ |c|c c c c  } 
 \hline
  & 0 & 1 & 2 & 3  \\ 
 \hline
 0 & $1$ & - & - & -  \\
 1 & - & $3$ & - & -   \\ 
 2 & - & - & $4$ &  $2$ \\ 
\end{tabular}
\end{center}
\end{thm}   

The proof of Theorem \ref{boocherhiyengarbetti} is also based on some more general results on the algebra structure on $H(R)$.  The authors of \cite[Theorem 3.3]{boocherhiyengar} prove that the diagonal subalgebra $\Delta(R)=\bigoplus_i H_i(R)_{2i}$ of $H(R)$ is a quotient of the exterior algebra on $H_1(R)_2$ by quadratic relations that depend only on the first syzygies of $I$.  As a consequence of this, they find that the linear strand and the main diagonal of the Betti table for a Koszul algebra satisfy the inequality in Question \ref{questbetti}, as we see in the following result.  We note that the proof of (2) also uses the characterization of complete intersections given by Bruns in Theorem \ref{bruns}. 

\begin{thm}(\cite[Corollary 3.4, Proposition 4.2]{boocherhiyengar})\label{bhici}
If $R$ is a Koszul algebra whose defining ideal is minimally generated by $g$ quadrics, then
\begin{description}[style=unboxed,leftmargin=0cm]
\item[(1)]  $\mathrm{dim}\,_kH_{i}(R)_{i+1}\leq{g\choose i}$ for $2\leq i\leq g$, and if equality holds for $i=2$, then I has height $1$ and a linear resolution of length $g$.
\item[(2)]  $\mathrm{dim}\,_kH_{i}(R)_{2i}\leq{g\choose i}$ for $2\leq i\leq g$, and if equality holds for some $i$, then I is a complete intersection.
\end{description}
\end{thm} 

We note that the inequality in Theorem \ref{bhici} (2) clearly holds for quadratic complete intersections, and the fact that it holds for Koszul algebras which are not complete intersections follows from the contrapositive of the statement in \cite[Corollary 3.4]{boocherhiyengar}.  In fact, this inequality also follows from the proofs of \cite[Theorem 3.1]{aci1} and \cite[Corollary 3.2]{aci1}.   

Building on the work of Boocher, Hassanzadeh, and Iyengar on the $g\leq 3$ case, Mantero and Mastroeni describe the Betti tables of Koszul algebras whose defining ideals are generated by $4$ quadrics in \cite{mastroeni}; their result is as follows.  

\begin{thm}(\cite[Main Theorem]{mastroeni})\label{mastroeni4}
Let $R=Q/I$ be a Koszul algebra with $I$ minimally generated by 4 quadrics and $\mathrm{ht}\,I=2$.  Then the Betti table of $R$ over $Q$ is one of the following: 
\begin{center}
\begin{tabular}{ |c|c c c c  } 
 \hline
  & 0 & 1 & 2 & 3  \\ 
 \hline
 0 & $1$ & - & - & -  \\
 1 & - & $4$ & $4$ & $1$   \\ 
\end{tabular}
\hspace{0.3cm}
\begin{tabular}{ |c|c c c c  } 
 \hline
  & 0 & 1 & 2 & 3  \\ 
 \hline
 0 & $1$ & - & - & -  \\
 1 & - & $4$ & $3$ & -   \\ 
 2 & - & - & $1$ &  $1$ \\ 
\end{tabular}
\hspace{0.3cm}
\begin{tabular}{ |c|c c c c c } 
 \hline
  & 0 & 1 & 2 & 3 & 4 \\ 
 \hline
 0 & $1$ & - & - & - & -  \\
 1 & - & $4$ & $3$ & $1$ & -  \\ 
 2 & - & - & $3$ &  $3$ & $1$ \\   
\end{tabular}
\hspace{0.3cm}
\begin{tabular}{ |c|c c c c c } 
 \hline
  & 0 & 1 & 2 & 3 & 4 \\ 
 \hline
 0 & $1$ & - & - & - & -  \\
 1 & - & $4$ & $2$ & - & -  \\ 
 2 & - & - & $4$ &  $4$ & $1$ \\   
\end{tabular}
\end{center}
\end{thm}

Theorem \ref{mastroeni4} gives an affirmative answer to \ref{questbetti} in the case that $g=4$ and $\mathrm{ht}\,I=2$; the case where $\mathrm{ht}\,I=4$ (i.e. $I$ is a complete intersection) is clear, the case where $\mathrm{ht}\,I=1$ follows from Theorem \ref{aci1mainthm}, and the case where $\mathrm{ht}\,I=3$ is settled by Mastroeni in \cite[Main Theorem]{mastroenialmost}.  Together these results give an affirmative answer to Question \ref{questbetti} when $g=4$.  However, Question \ref{questbetti} remains open for $g>4$. 

\subsection{Subadditivity of syzygies}

The focus of this section is the following conjecture of Avramov, Conca, and Iyengar from \cite{aci2} concerning the subadditivity of syzygies, and equivalently the subadditivity of Koszul homologies.

\begin{conj}(\cite[Conjecture 6.4]{aci2})\label{conj}
If $R$ is a Koszul algebra, then the following inequality holds
\begin{align*}
t_{a+b}(R)\leq t_a(R)+t_b(R)
\end{align*}
whenever $a+b\leq\mathrm{pd}_QR$.
\end{conj}  

This conjecture is based on the following result, which asserts that a weaker inequality holds.

\begin{thm}(\cite[Theorem 6.2]{aci2})\label{subadd}
Let $a$ and $b$ be non-negative integers such that ${a+b}\choose a$ is invertible in $k$ and $\mathrm{max}\{a,b\}\leq m(R)$.  If $R$ is Koszul, then the following inequalities hold
\begin{align*}
t_{a+1}(R)&\leq t_a(R)+2 \\
t_{a+b}(R)&\leq t_a(R)+t_b(R)+1
\end{align*}
for $b\geq 2$.
\end{thm}

In fact, both Conjecture \ref{conj} and Theorem \ref{subadd} are stated more generally in \cite{aci2}; rather than requiring that $R$ is Koszul, it is only required that the minimal resolution of $k$ is linear up to a certain degree.

To our knowledge, Conjecture \ref{conj} remains open; however, the following results of Eisenbud, Huneke, and Ulrich and of Herzog and Srinivasan provide supporting evidence for its validity.

\begin{thm}\label{evidence}
Conjecture \ref{conj} holds in the following cases:
\begin{description}[style=unboxed,leftmargin=0cm]
\item[(1)]  (\cite[Corollary 4.1]{ehu}) $\mathrm{dim}R\leq 1$, $\mathrm{depth}R=0$, and $a+b=\mathrm{rank}_kR_1$; 
\item[(2)]  (\cite[Corollary 2.1]{FRG},\,\cite[Corollary 4]{HS}) $R$ has monomial relations and $b=1$;  
\end{description}
\end{thm}

We note that the quadratic monomial case in Theorem \ref{evidence} (2) is treated in \cite{FRG} and the general monomial case is treated in \cite{HS}.  

Furthermore, in \cite{HS} Herzog and Srinivasan prove the following case of Conjecture \ref{conj} (without the Koszul assumption), which improves a result of McCullough in \cite{mccullough}.  

\begin{thm}(\cite[Corollary 3]{HS},\,\cite[Theorem 4.4]{mccullough})
Let $\mathrm{pd}_Q\,R=p$. Then the following inequality holds
\begin{align*}
t_p(R)\leq t_{p-1}(R)+t_1(R).
\end{align*} 
\end{thm}

We finish this section by noting that the question of when the subadditivity property holds is an interesting one, even for non-Koszul algebras.  The subadditivity property holds for complete intersection ideals, but not for Gorenstein ideals; these are recent results of McCullough and Seceleanu in \cite{ms}.  Although the subadditivity property holds for certain classes of monomial ideals (see for example \cite{faridi}), the question remains open for general monomial ideals (as well as for Koszul algebras, as discussed above).

\subsection{Towards a characterization}

Although the Koszul property of $R$ has not been shown to be characterized completely by the algebra structure on $H(R)$, recent work of Myers in \cite{myers} introduces a property of $H(R)$ that is equivalent to the Koszul property of  $R$ under some extra conditions.  Myers defines $H(R)$ to be \textit{strand-Koszul} if $H'(R)=\bigoplus_{j-i=n} H_i(R)_j$ is Koszul in the classical sense; that is, if $k$ has a linear resolution over $H'(R)$.  With this terminology, he proves the following result.

\begin{thm}(\cite[Theorem C]{myers})\label{myers}
Let $R$ be a standard graded $k$-algebra with $k$ a field.  If $H(R)$ is strand-Koszul, then $R$ is Koszul.  Furthermore, the converse holds if $R$ satisfies one of the following conditions:
\begin{description}[style=unboxed,leftmargin=0cm]
\item[(1)] $R$ has embedding dimension $\leq 3$;
\item[(2)] $R$ is Golod;
\item[(3)] The defining ideal of $R$ is minimally generated by 3 elements;
\item[(4)] $R$ is a quadratic complete intersection;
\item[(5)] The defining ideal of $R$ is the edge ideal of a path on $\geq 3$ vertices.
\end{description}
\end{thm} 

Despite the connections we have seen between the Koszul property and the algebra structure on Koszul homology throughout this section and the progress towards a characterization highlighted in the theorem above, results of Roos in \cite{roos} present potential challenges, or perhaps even obstructions, to finding a characterization of Koszul algebras by the algebra structures on their Koszul homologies alone, as we see in the following theorem.  

\begin{thm}(\cite[Section 3]{roos})\label{thmroos}
The Koszul homology algebras of the Koszul algebras 
\begin{align*}
R_{71}&=k[x,y,z,u]/(x^2,\, y^2,\, z^2,\, u^2,\, xy,\, zu,\, yz-xu) \\ 
\end{align*}
and 
\begin{align*}
R_{71v16}&=k[x,y,z,u]/(xz+u^2,\, xy,\, xu,\, x^2,\, y^2+z^2,\, zu,\, yz),
\end{align*}
are isomorphic as bigraded $k$-vector spaces; however, they are not isomorphic as $k$-algebras.
\end{thm}

Roos uses the DGAlgebras package written by Moore to describe and compare the Koszul homology algebra structures of the Koszul algebras $A=R_{71}$ and $B=R_{71v16}$ defined above.  Roos shows that the Poincar{\'e} series $P_k^A(t)$ and $P_k^B(t)$, and the Betti tables of $A$ and $B$ over $Q=k[x,y,z,u]$ are both the same, hence $H(A)$ and $H(B)$ are isomorphic as bigraded $k$-vector spaces; however, the algebra structures on $H(A)$ and $H(B)$ are different.  In fact, Roos concludes that the Eilenberg-Moore spectral sequence degenerates over $H(A)$ but not over $H(B)$; thus, by \cite[Theorem B]{myers}, $H(A)$ is strand-Koszul and $H(B)$ is not.

We end by noting that the algebra $B=R_{71v16}$ provides an example of a Koszul algebra whose Koszul homology algebra is not strand-Koszul; thus, the converse of Theorem \ref{myers} does not hold in general.    
  

\bibliographystyle{siam}
\bibliography{survey}
\nocite{*}

\vspace{1.25cm}

\noindent \textsc{Department of Mathematics, Yale University, New Haven, CT 06520}
\vspace{2pt}
\newline
\textit{E-mail Address}: rachel.diethorn@yale.edu

 \end{document}